\documentclass{article}

\usepackage[latin1]{inputenc}
\usepackage[square,numbers]{natbib}
\usepackage[english]{babel}

\makeindex

\everymath{\displaystyle} 
 
\usepackage{calc}
 
\usepackage{lmodern} 
\usepackage{theorem,latexsym,amsmath,amsfonts,amssymb,mathrsfs}
\usepackage[T1]{fontenc}
\usepackage{fancybox}
\usepackage{fancyhdr} 
\usepackage[english]{babel} 
\usepackage{multicol,multirow}
\usepackage{float}
\usepackage{array}
\usepackage[normalem]{ulem}
\usepackage{color}
\usepackage{tablists}
\usepackage{subfloat}
\usepackage{ragged2e}
\usepackage{graphicx}
\usepackage{geometry}
\usepackage{indentfirst}
\usepackage{latexsym}

\usepackage{sectsty} 
\usepackage{titlesec} 
\usepackage{enumerate} 
\usepackage{cancel} 
\usepackage{subfloat} 
\usepackage{sectsty} 
\usepackage[colorlinks, citecolor=blue,linkcolor=black]{hyperref}

\titleformat*{\subsection}{\bfseries\rmfamily}
\titleformat*{\section}{\bfseries\rmfamily}

\theorembodyfont{\slshape}

\usepackage[para]{footmisc}

\titlelabel{\thetitle.~}

\numberwithin{equation}{section}

\AtBeginDocument{}

\newtheorem{Lemma}{Lemma}[section]{\bfseries}{\itshape}
\newtheorem{Proposition}{Proposition}[section]{\bfseries}{\itshape}

\newtheorem{Remark}{Remark}

{\bfseries}{\itshape}
\newtheorem{Assumption}{Assumption} 
\newtheorem{Theorem}{Theorem}

\newcommand{\R}{\mathbb R}

\newcommand{\N}{\mathbb N}
\newcommand{\p}{\mathbb P}
\newcommand{\E}{\mathbb E}

\title{
\textbf{Asymptotic normality and strong consistency of kernel regression estimation in q-calculus}}
\date { }
\author{\textbf{Emmanuel De Dieu Nkou}${}^{(a)}$ and \textbf{Fridolin Melong}${}^{(b)}$\\
\\
$\mathbf{{}^{(a)}}$\textit{Laboratoire de Probabilit\'es, Statistique et Informatique (LPSI),}\\
 \textit{Unit\'e de Recherche en Math\'ematiques et Informatique (URMI),}\\
\textit{Universit\'e des Sciences et Techniques de Masuku, BP 813 Franceville, Gabon}\\
mail: \textit{emmanueldedieunkou@gmail.com}
\\
\\
$\mathbf{{}^{(b)}}$\textit{International Chair in Mathematical Physics and Applications (ICMPA-UNESCO Chair),}\\
\textit{University of Abomey-Calavi, 072 BP 50 Cotonou, Republic of Benin}\\
mail: \textit{fridomelong@gmail.com}
}

\begin{document}

\maketitle

\medskip

\begin{center}
\textbf{Abstract}
\end{center}
We construct a family of estimators for a regression function based on a sample following a $q$-distribution. Our approach is nonparametric, using kernel methods built from operations that leverage the properties of $q$-calculus. Furthermore, under appropriate assumptions, we establish the weak convergence and strong consistency of this family of estimators.

\medskip

\begin{center}
\textbf{R\'{e}sum\'{e}}
\end{center}
Nous construisons une famille d'estimateurs pour une fonction de r\'{e}gression bas\'{e}e sur un \'{e}chantillon suivant une $q$-distribution . Notre approche est non param\'{e}trique, utilisant des m\'{e}thodes \`{a} noyau construites \`{a} partir d'op\'{e}rations qui exploitent les propri\'{e}t\'{e}s du $q$-calcul. De plus, sous des hypoth\`{e}  ses appropri\'{e}es, nous \'{e}tablissons la convergence faible et la convergence presque s\^{u}re de cette famille d'estimateurs.

\bigskip

\textbf{Keywords:} $q$-calculus, nonparametric estimation, regression estimation.

\textbf{2020Mathematics Subject Classification}: 05A30, 62G05, 62G07.

\vspace{1cm}

\section{Introduction}

The quantum calculus (or $q$-calculus) emerged as a link between mathematics and physics. It turns out to be even more significant than classical computation. It depends on a parameter $q>0$. It encompasses classical computation: when $q\rightarrow 1$, the various results obtained reduce to classical operations. The $q$-calculus can thus be considered as an extension of the classical calculus. A another example can be found in classical analysis, where the definition of the derivative relies on the existence of a limit, in contrast, in $q$-calculus, the notion of $q$-derivatives can exist without requiring the definition of a limit. It has many applications in different areas of mathematics such as number theory, combinatorics, orthogonal polynomials, basic hypergeometric functions, quantum theory and electronics, and in different areas of applied mathematics such as probability and statistics in which it abounds \citep{bensalah2021, bouzida2020, hilhorst2010,  umarov2008}. In the specific case of probability and statistics, two types of functions are the focus of a great deal of research due to their importance: density and regression functions. Their significance arises from their applications in many areas of probability, statistics, and even beyond \citep{ebende2024, zhu1996, zhu2007}. It is therefore clear that estimating the density and regression function using various methods is a crucial task in the field of statistics and probability. However, despite the abundance of statistical results related to $q$-calculus, there is very few which address to nonparametric statistics using kernel methods in the estimation of density functions and regression functions. Nonparametric estimation by kernel method using $q$-calculus is performed in \cite{badrani2024}. The authors, with $q$-calculus methods and in the case where $0<q<1$, obtain   the $q$-asymptotic normality of kernel density estimator.In this paper, we adopt a similar approach to the study of the asymptotic normality of the density in \cite{badrani2024} and extend this work. More precisely, under relaxed assumptions, we investigate the strong convergence of this estimator and determine its rate of convergence. Furthermore, we apply the same $q$-calculus operations to investigate the asymptotic normality and almost sure convergence of the regression function estimation using the kernel method, referred to as the $q$-kernel method.

Let the pair $(X,\, Y)$ be a pair of random variables which follow a $q$-distribution in $q$-probability space $\left(\Omega, \, \mathcal{F}, \p_q\right)$.  We denote by $r$ the regression function of the random variable $Y$ on the random variable $X$ defined by
$$
r(x)=\E_q\left(Y|X=x\right).
$$
$\E_q(\cdot)$ is $q$-analog expectation in $q$-calculus. It will be defined in the next section. We approach the estimation of this function using a sample of independent and identically distributed random variables that follow a $q$-distribution similar to that of the pair $(X,Y)$. The estimation is performed without making any prior assumptions about the distribution of $(X,Y)$, except on the existence of the regression function $r$. The estimate taken here is that of Nadaraya-Watson. To achieve optimal results, we construct a family of $q$-kernels that converge to a family of classical kernels when $q$ approaches 1

This paper is organized as follows: in Section 2, we introduce some concepts from quantum calculus that are useful for this study, along with essential results that support our work. In Section 3, we present $q$-kernel estimators for $f$ and $r$, as well as some auxiliary functions. In Section 4, we outline the assumptions used and state the main results of this work. Finally, Section 5 is dedicated to proving complementary results employed in our analysis.

\section{Preliminaries}

In this section, we briefly recall main definitions, notations of  some basics of $q$-calculus that will be useful to us in this paper. These notions are presented in \citep{badrani2024, boutouria2018, diaz2009} among many others.

Throughout this paper we suppose $0 < q < 1$. The $q$-analog $[n]_q$ of any positive integer $n\in \N^*$ and $q$-analog factorial $[n]_q!$  are defined as 
$$
[n]_q=\frac{1-q^n}{1-q}\, \mbox{ and }\, [n]_q!=\left[n\right]_q \times \left[n-1\right]_q\times\cdots \times [1]_q.
$$
The $q$-analog of $(x-a)^n$ is
$$
(x-a)_q^n=
\left\{
\begin{array}{lll}
1 && \mbox{if } n = 0 \\
(x-a) \times (x-qa)\times\cdots \times (x-q^{n-1}a)&& \mbox{if } n\geqslant 1,
\end{array}
\right.
$$
and we have
$$
(a-x)_q^n=(-1)^nq^{\frac{n(n-1)}{2}}\left(x-q^{1-n}a\right)_q^n.
$$
The $q$-analog of the exponential function $e^x$ given by
$$
e_q^x = \sum_{k=0}^{+\infty}\frac{x^k}{[k]_q!}.
$$
The $q$-analog of identity $e^x e^{-x}=1$ is defined by $e_q^x E_q^{-x}=1$, where
$$
E_q^{x}= e_{\frac{1}{q}}^x = \sum_{k=0}^{+\infty} q^{\frac{k(k-1)}{2}}\frac{x^k}{[k]_q!}.
$$
The $q$-analog of $e^{-\frac{x^2}{2}}$ is given by
$$
E_{q^2}^{-\frac{q^2x^2}{[2]_q}} = \sum_{k=0}^{+\infty}\frac{q^{k(k+1)}(q-1)^k}{(1-q^2)_{q^2}^k}   x^{2k}.
$$
The $q$-analog of an improper integral is proper integral with limits $-\nu$ and $\nu$ where
$$
\nu=\nu(q)=\frac{1}{\sqrt{1-q}}.
$$
For $a,\,b\in \R$, the  Jackson integral or $q$-integral of arbitrary function $f:\R \rightarrow \R$ on $[a,b]$ is defined by
\begin{equation}\label{eqnint1}
\int_a^b f(x)d_qx = (1-q)\sum_{k=0}^{+\infty}q^k \left[bf(q^kb) - a f(q^ka)\right],
\end{equation}
When $f$ is continue on $[a,b]$, we get that
$$
\lim_{q\rightarrow 1} \int_a^b f(x)d_qx = \int_a^b f(x)dx.
$$
We have the $q$-derivative of a function $f$ by the following expression
$$
D_qf(x)=\frac{d_qf(x)}{d_qx}=\frac{f(qx) - f(x)}{qx-x}.
$$
Clearly, if $f$ is differentiable, then
$$
\lim_{q\rightarrow1}D_qf(x)=f'(x).
$$
Then, we can define the $s$\, $q$-derivatives of $f$ by
$$
D_q^sf(x) =
\left\{
\begin{array}{lll}
f(x) && \mbox{if } s = 0 \\
D_q\left(D_q^{s-1}f\right)(x)&& \mbox{if } s\geqslant 1.
\end{array}
\right.
$$
To derive some of our results, we will use the $q$-Taylor formula for a function $f$ being differentiable $s$ times, which can be found in \cite{ernst1999, marinkovic2002}. This  formula  is given by
\begin{equation}\label{eq1taylor}
f(b)= \sum_{k=0}^s \frac{(b-a)_q^k}{[k]_q!}\left(D_q^k f\right)(a) + R_s(f,b,a,q),
\end{equation}
where $R_s(f,b,a,q)$ is the remainder term determined by
$$
R_s(f,b,a,q) = \int_a^b \frac{(b-t)_q^s}{b-t}\frac{\left(D_q^k f\right)(t)}{[s-1]_q!}d_qt.
$$

\bigskip

Given a random variable $X$, one defines the positive and negative parts by $X_{+} = \max(X, 0)$ and $X_{-} =\min(X, 0)$. These are non negative random variables, and it can be directly checked that $X=X_{+}+X_{-}$.

The $q$-logarithm and $q$-exponential functions were originally defined in \cite{tsallis1994} by
$$
\ln_q(x)=\frac{x^{1-q}-1}{1-q}\,\mbox{ and }\,\exp_q(x)=\left[1+(1-q)x\right]^{\frac{1}{1-q}}\,\mbox{ if }\,1+(1-q)x \geqslant 0,
$$
with $\lim_{q\rightarrow 1}\ln_q(x) = \ln(x):=\ln_1(x)$ and $\lim_{q\rightarrow 1}\exp_q(x) = \exp(x):=\exp_1(x)$.

The following result will be considered:
\begin{equation}\label{eqnlnq1}
\mbox{ For } \, q\in [0,1],\quad \ln_q(x+1)\leqslant x.
\end{equation}

This result is clear for $q=1$ or for $q=0$. Now, for $q\in]0;1[$, define for $x\in ]-1;+\infty[$ the function $\tau_q(x) = \ln_q(1+x) -x=\frac{(1+x)^{1-q}-1 - x(1-q)}{1-q}$. \\
We have $\frac{d\tau_q(x)}{dx} =(1-q)\left((1+x)^{-q} - 1\right)$. A simple study of variation shows us $\sup_{x\in ]-1;+\infty[} \tau_q(x)=0$, then $\tau_q(x)\leqslant 0$.

Note that if $xy=qyx$, we have 
\begin{equation}\label{eqnlnq2}
\ln_q(xy)= \ln_q(x)+ \ln_q(y)\, \mbox{ and }\, \exp_q(x)\exp_q(y)=\exp_q(x+y).
\end{equation}
From (\ref{eqnlnq2}) we can deduce by induction that for all integer $n\in \N$
\begin{equation}\label{eqnlnq3}
\ln_q(x^n)= n\ln_q(x)\, \mbox{ ,  }\, \left(\exp_q(x)\right)^n =\exp_q(nx) \, \mbox{ and  }\,\left(\exp_q(x)\right)^{-n} =\exp_q(-nx).
\end{equation}

\medskip
To conclude this section, we introduce the space $L_q^{\infty}\left(\R\right)$, which will contain two essential functions necessary for obtaining our results. First we present the space $L_q^{p}\left(\R\right)$ for $1 \leqslant p<\infty$. The space $L_q^{p}\left(\R\right)$ is defined as follows
$$
L_q^{p}\left(\R\right) = \left\{f:\R \rightarrow \R\,\mbox{ measurable such that }\, \int_\R |f(x)|^p d_qx< \infty\right\}.
$$
Clearly, if $f\in L_q^{p}\left(\R\right)$, the associated $q$-norm is expressed by  $f$ is $\|f\|_{L_q^{p}} = \left(\int_\R |f(x)|^p d_qx\right)^{\frac{1}{p}}$.\\
From this fact, the space $L_q^{\infty}\left(\R\right)$ is defined as follows
$$
L_q^{\infty}\left(\R\right) = \left\{f:\R \rightarrow \R\,\mbox{ measurable such that }\, \|f\|_{L_q^{\infty}}=\sup_{x\in \R}|f(x)|< \infty\right\}.
$$

\bigskip

\section{The $q$-kernel regression estimation}
Let $\left(X,\, Y\right)$ be a two dimensional random variables on $q$-probability space $\left(\Omega, \mathcal{F},\p_q\right)$. The $q$-regression function of $Y$ given $X$, assuming that $\E_q\left(|Y|\right)<\infty$, is defined by 
$$
r(x)=\E_q\left(Y | X=x\right) = \frac{g(x)}{f(x)},
$$
 $g$ is function defined by $g(x) = r(x)f(x)=\int_{\R}yf_{_{(X,Y)}}(x,y)d_qy$ where  $f_{_{(X,Y)}}$ is $q$-density function of pair  $\left(X,\, Y\right)$. $f$ is $q$-density function of $X$, which means that 
$$
\forall x \in \R,\,  f(x)> 0, \quad \mbox{ and }\quad \int_{\R} f(x)d_qx=1.
$$
The $q$-expectation of a random variable $X$ is presented as follows
$$
\E_q(X)= \int_{\R}x f(x)d_qx,
$$
and the $q$-variance of  $X$ as follows
$$
Var_q(X)= \E_q\left[\left(X-\E_q(X)\right)^2\right]= \E_q\left(X^2\right) - \left(\E_q(X)\right)^2.
$$
Note that the $q$-expected value operator $\E_q(\cdot)$ is linear in the sense that
$$
\E_q(X_1 + X_2)=\E_q(X_1) + \E_q(X_2),
$$
and if $X_1$ and $X_2$ are independent random variables, then
$$
\E_q(X_1 \, X_2)=\E_q(X_1) \, \E_q(X_2).
$$

\begin{Remark}
It is well known that, in the construction of $r$, $X$ and $Y$ are not independent. Therefore, knowing the value of $X$ influences the realization of $Y$, typically leading to a decrease in its value. $r(X)$ is the random variable that is as close as possible to $Y$ in quadratic mean. That is to say, it is the one that minimizes $\E_q\left[\left(Y - r(X)\right)^2\right]$.
\end{Remark}

\medskip

\begin{Remark}
In the classic case of regression function $r(x)$, due to his quotient form, controlling its denominator for very small values is essential to obtaining valid results. For overcoming this technical difficulty, some studies have introduced assumptions on $f(x)$, such that $f(x)>c$ with $c>0$ \citep{zhu2007} or by replacing $f(x)$ with $f_{b_n}(x)=\max\{b_n,f(x)\}$ where $b_n$ is a sequence  tending to $0$ as $n\rightarrow + \infty$ \citep{zhu1996, nkou2019}. It is clear that the second condition provides an improvement over the first. In this paper, we propose a condition that unifies these two approaches, namely:
\begin{equation}\label{eqnbnf1}
f(x)\geqslant b_n \, \mbox{ where } \,  b_n \rightarrow 0\,\mbox{ as } \, n\rightarrow +\infty.
\end{equation}
We will simply consider this condition in this paper.
\end{Remark}

We will subsequently give an estimator of $r$ by the kernel method. Let $\left\{\left(X_1,Y_1\right),\, \cdots, \left(X_n,Y_n\right)\right\}$ be $n$ independent and identically $q$-distributed of $\left(X,\, Y\right)$ on the common $q$-probability space $\left(\Omega, \mathcal{F},\p_q\right)$. The kernel density estimator  for unknown $q$-regression $r$ at any point $x$ is defined as 
$$
\widehat{r}_n(x)=\frac{\widehat{g}_n(x)}{\widehat{f}_n(x)},
$$
where
$$
\widehat{f}_n(x) = \frac{1}{nh_n}\sum_{i=1}^n K_q\left(\frac{x-X_i}{h_n}\right)\, \mbox{ and }\, \widehat{g}_n(x) = \frac{1}{nh_n}\sum_{i=1}^n Y_iK_q\left(\frac{x-X_i}{h_n}\right).
$$
$h_n$ is a numerical sequence, as in the classical case, it must first satisfy the condition $h_n \rightarrow 0$ as $n\rightarrow +\infty$. Starting from this condition, we introduce additional conditions in the next section that are necessary to obtain our results. $K_q(x)$  is a functional kernel that depends on the parameter $q$. Due to this dependency, we refer to it as the $q$-kernel. To ensure that the $q$-kernel function $K_q(\cdot)$ satisfies certain properties that are essential to obtain our results,  it could be a  $q$-Gaussian kernel
$$
K_q(u)=\frac{1}{c(q)}E_{q^2}^{-\frac{q^2u^2}{[2]_q}},\quad u\in[-\nu,\nu], 
$$
where $c(q)=2(1-q)\nu \sum_{k=0}^{+\infty}q^k E_{q^2}^{-\frac{q^2(q^k\nu)^2}{[2]_q}}$ is $q$-analog $\sqrt{2\pi}$.

$K_q(\cdot)$ could also be one the functions belongs to the family of kernels $K_{q,p}(x)$  in the paragraph below.

We propose, starting from the classical kernels, some versions of $q$-kernels and their different analytical properties. 
In $[-1,1]$, following (\ref{eqnint1}), the Jackson integral or $q$-integral, of an  function $\gamma:\R \rightarrow \R$ is defined by
\begin{equation}\label{eqn1qint}
\int_{-1}^1 \gamma(x)d_qx = (1-q)\sum_{k=0}^{+\infty} q^k\left[\gamma\left(q^k\right) + \gamma\left(-q^k\right)\right].
\end{equation}
Let consider the family of functions $\gamma_{q,p}(x)=\left(1-q^2x^2\right)^p 1_{\{|x|\leqslant 1\}}$, where $p\in \N^*$ and the family of kernels defined by
$$
K_{q,p}(x)=\frac{\gamma_{q,p}(x)}{\int_{-1}^1\gamma_{q,p}(x)d_qx},\quad 0<q<1.
$$
We have $\gamma_{q,p}\left(q^k\right)=\gamma_{q,p}\left(-q^k\right) = \sum_{\ell=0}^p (-1)^\ell \binom{p}{\ell}q^{2(k+1)\ell}$. Then, by (\ref{eqn1qint}), we have
$$
\int_{-1}^1\gamma_{q,p}(x)d_qx =2(1-q) \sum_{k=0}^{+\infty} q^k\sum_{\ell = 0}^p (-1)^\ell \binom{p}{\ell}q^{2(k+1)\ell} = 2 \sum_{\ell = 0}^p (-1)^\ell \binom{p}{\ell} \frac{q^{2\ell}}{[2\ell+1]_q}.
$$
Consequently
\begin{equation}\label{eqn1kqp}
K_{q,p}(x)=\frac{1}{c_q}\gamma_{q,p}(x),\mbox{ where }\, c_q = 2 \sum_{\ell = 0}^p (-1)^\ell \binom{p}{\ell} \frac{q^{2\ell}}{[2\ell+1]_q},\quad 0<q<1.
\end{equation}
We also see that, by using (\ref{eqn1qint}) on the function $\gamma(x)=x\, \gamma_{q,p}^m(x)$, with $m\in \N^*$, we will have $\gamma\left(q^k\right) + \gamma\left(-q^k\right)=0$, which means that \, $\int_{-1}^1 x\, \gamma_{q,p}^m(x)d_qx =0$. And therefore
$$
\int_{-1}^1 x K_{q,p}^m(x) d_qx=0.
$$
Some simple examples of $q$-kernel of family $K_{q,p}$:
\begin{description}
	\item[$\bullet$ $q$-rectangular kernel ($p=0$):] $K_q(x)=\frac{1}{2}1_{\{|x| \leqslant 1\}}$.
The $q$-analog rectangular kernel is the same that the classic rectangular kernel.
\item[$\bullet$ $q$-Epanechnikov kernel ($p=1$):] $K_q(x)=\frac{1}{c_q}\left(1-q^2x^2\right)1_{\{|x| \leqslant 1\}}$, $c_q=\frac{2([3]_q-q^2)}{[3]_q}$;
\item[$\bullet$ $q$-Biweight kernel ($p=2$):]  $K_q(x)=\frac{1}{c_q}\left(1-q^2x^2\right)^2 1_{\{|x| \leqslant 1\}}$, $c_q=2\left(1-\frac{2q^2}{[3]_q}+\frac{q^4}{[5]_q}\right)$.
\end{description}
Note that $\lim_{q\rightarrow 1} K_q(x)= K_1(x)$ is verified. $K_1(x)$ is a classic kernel. For $p = 0,\, 1, \,2, \,3$, $K_1(x)$ stands one of  resulting kernels are known as the uniform, Epanechnikov (or quadratic), biweight (or quartic) and triweight kernels respectively.

\bigskip

\section{Assumptions and results}

In this section, we first present the assumptions, followed by the results establishing the asymptotic normality of the $q$-kernel estimator $\widehat{r}_n (x)$. We then conclude with the results demonstrating the almost sure convergence of $\widehat{f}_n (x)$ and $\widehat{r}_n (x)$. Additional results regarding the bias and variance are presented  in Section \ref{proofs1}.

\subsection{Assumptions}
\begin{Assumption}\label{hyph1}
The sequences  $h_n$ and $b_n$ are such that:
 $$h_n b_n^{-2} \rightarrow 0\, \mbox{ , } \, nb_n^2 h_n \rightarrow +\infty \, \mbox{ , }\, \frac{\ln_q n}{n b_n h_n}  \rightarrow 0\, \mbox{ as } \, n\rightarrow + \infty.$$
\end{Assumption}
\begin{Assumption}\label{hypf1}
The functions $f$  and $g$ belong to $L_q^{\infty}\left(\R\right)$ and are  three times $q$-differentiable and their respective first, second and third differentials check, for all $x\in\R$
$$
\max_{1\leqslant k \leqslant 3}\left(\left|(D_q^{(k)} f)(x)\right| \, ,  \,\left|(D_q^{(k)} g)(x)\right| \right) < +\infty.
$$
\end{Assumption}

\begin{Assumption}\label{hypy1}
There exists a sequence $A_{nq}$ of strictly positive numbers such that $A_{nq} \sim \sqrt{\ln_q (n)}$ and $\max_{1\leqslant i \leqslant n} |Y_i| \leqslant A_{nq}$.

\end{Assumption}
\begin{Assumption}\label{hypk1}
$K_q$ is the $q$-kernel function, which is assumed to satisfy the following properties:
\begin{enumerate}[1.]
	\item $\sup_{u\in\R} K_q(u) \leqslant M$,\, $\int_{-\nu}^\nu K_q(u)d_qu=1$ and $\int_{-\nu}^\nu u K_q^m(u)d_qu=0$ for $m\in \N^*$;
	\item $\int_{-\nu}^\nu u^2 K_q(u)d_qu <+\infty$ and $\int_{-\nu}^\nu  K_q^3(u)d_qu <+\infty$.
\end{enumerate}
\end{Assumption}

\begin{Remark}
In Assumption \ref{hyph1}, the conditions on $h_n$ and $b_n$ are more general compared to those in the classical case (see \cite{zhu2007}). A future study may be conducted to determine the optimal choice for this sequence. The sequence $b_n$ is defined as defined as $b_n = \max(\varepsilon,\, n^{-c})$  to obtain an accurate estimate of $\widehat{f}_n$ and $\widehat{r}_n$; more details can be found in \cite{nkou2019}. In Assumption \ref{hypf1}, we provide further details on these commonly accepted bounds, which are essential for obtaining results. Assumption \ref{hypk1} is weaker than the boundedness assumption; for instance, it has been considered in \cite{nkou2023}. In Assumption \ref{hypk1}, to ensure that the $q$-kernel function $K_q(\cdot)$ satisfies these properties mentioned, we propose the $q$-Gaussian kernel and the family of $q$-kernels $K_{q,p}(\cdot)$ defined in (\ref{eqn1kqp}). 
\end{Remark}

\subsection{Main results}

In the first theorem,  we establish the asymptotic normality of $\widehat{r}_n(x) - r(x)$. An estimator is asymptotically normal in $q$-calculus meaning that its asymptotic distribution is the normal $q$-distribution.It is a consequence of Proposition \ref{propbr1} and \ref{propo2an}, Theorem \ref{theoly1} and  Lemma \ref{lemaz1}.

\begin{Theorem}\label{Theonr1} Under Assumptions \ref{hyph1}, \ref{hypf1}, \ref{hypy1} and \ref{hypk1}, if $\sqrt{nh_n^5}\rightarrow c$
$$
\sqrt{nh_n}\left(\widehat{r}_n (x) - r(x)\right) \stackrel{\mathscr{L}}{\longrightarrow } \mathcal{N}_q\Big(\mathscr{E}_q(x)\, , \, \mathscr{V}_q(x)\big)\,\mbox{ as }\, n\rightarrow +\infty,
$$
where
\begin{eqnarray*}
\mathscr{E}_q(x) &=& \frac{q c}{[2]_q f (x)} \Big(\left(D_q^2 g\right)(x) - r(x)\left(D_q^2 f\right)(x)\Big)\int_{-\nu}^{+\nu} u^2K_q(u)d_qu\\
\mathscr{V}_q(x) &=& f(x)^{-1}\left(Var_q \left(Y | X=x\right)\int_{-\nu}^{+\nu} K_q^2(u)d_qu\right).
\end{eqnarray*}
\end{Theorem}

\medskip

The following theorem presents, while highlighting the rate of convergence, the strong consistency  of the estimators $\widehat{f}_n$ and $\widehat{g}_n$.\\
We define for $k\in \N$
\begin{equation}\label{eqnrho1}
\Gamma_k(x)=\int_{\R}y^k f_{_{(X_t,Y_t)}}(x,y)d_qy = \E_q\left(Y^k | X=x\right)f(x).
\end{equation}
Clearly, $\Gamma_0(x)=f(x)$ and $\Gamma_1(x)=g(x)$, which implies that the $q$-kernel estimator of $\Gamma_k$ with the sequence of random variables $\left\{\left(X_1,Y_1\right), \, \left(X_2,Y_2\right),\, \cdots, \left(X_n,Y_n\right)\right\}$ is 
$$
\widehat{\Gamma}_{k,n}(x)=\frac{1}{nh_n}\sum_{i=1}^n Y_i^k K_q\left(\frac{x-X_i}{h_n}\right).
$$
We can therefore clearly see that $\widehat{\Gamma}_{0,n}(x) = \widehat{f}_n$ is the $q$-kernel density estimator of $f$ and $\widehat{\Gamma}_{1,n}(x)=\widehat{g}_n$ is $q$-kernel estimator of $g$.

We assume that, with Assumption \ref{hypy1} above, there exists $c_0>0$ such that, for $n$ large enough and $k\in \N$ 
$$
\|\Gamma_k\|_{L_q^{\infty}}=\sup_{x\in\R}|\Gamma_k(x)|< \big(c_0 \ln_q (n)\big)^{\frac{k}{2}}\|f\|_{L_q^{\infty}}.
$$

\medskip

Before presenting the next results, the following notations are necessary
\begin{equation}\label{eqn1vnk}
\begin{split}
v_{(n,k)}& = \frac{1}{nh_n}\left(\left(c_0  \ln_q (n)\right)^{\frac{k}{2}}\|f\|_{L_q^{\infty}} + o(1)\right) \int_{-\nu}^{+\nu} K_q^2(u)d_qu \\
w & = \frac{L}{[2]_q} \int_{-\nu}^{+\nu} u^2K_q(u)du + o(1). \\
\end{split}
\end{equation}

\begin{Theorem}\label{theoconstg1} Under Assumptions  \ref{hyph1}, \ref{hypf1}, \ref{hypy1} and  \ref{hypk1}, for all $t>0$, for $k\in \{0,1\}$ and $n\geqslant 1$, we have
$$
\sup_{x\in\R}\left|\widehat{\Gamma}_{k,n}(x)-\Gamma_{k}(x)\right|= O_{a.s}\left(qw h_n^2 + \sqrt{v_{(n,k)}\ln_q n}\right).
$$
\end{Theorem}
We refer to (\ref{eqn1eqf}) and Proposition \ref{propqg1} for proof of this Theorem.

\medskip

From the result of the previous theorem, we can easily obtain the strong convergence of the estimator $\widehat{r}_n$. This result is illustrated in Theorem \ref{theoconstr1} below.

\begin{Theorem}\label{theoconstr1} Under Assumptions  \ref{hyph1}, \ref{hypf1}, \ref{hypy1} and  \ref{hypk1}, for all $t>0$ and $n\geqslant 1$, we have
$$
\sup_{x\in\R}\left|\widehat{r}_n(x)-r(x)\right|= O_{a.s}\left(qw h_n^2 + \sqrt{v_{(n,1)}\ln_q n}\right).
$$
\end{Theorem}

\medskip

\subsection{Auxiliary results}

To obtain the results presented above, we used auxiliary results, which are outlined in this section. This includes the determination of biases and variances of the used estimators. To achieve this, we express $\widehat{r}_n(x) - r(x)$ by using another known random variable that we will denote by $Z_i(x)$. It is defined by
$$
Z_i(x) = \frac{1}{nh_n}\big(Y_i - r(x)\big) K_q\left(\frac{x-X_i}{h_n}\right).
$$
From this random variable, we derive the first results. Clearly, for all $x\in \R$ such that (\ref{eqnbnf1}) is hold, we have the following equality:
$$
\widehat{r}_n (x) - r(x) = \frac{\mathcal{A}_n(x)}{\widehat{f}_n (x)},
$$
where
\begin{equation}\label{eqn1anz}
\mathcal{A}_n(x)= \frac{1}{n}\sum_{i=1}^n\left\{\big(Y_i - r(x)\big)\frac{1}{h_n} K_q\left(\frac{x-X_i}{h_n}\right)\right\}=\sum_{i=1}^n Z_i(x). 
\end{equation}
Then, the asymptotic behavior of $\widehat{r}_n (x) - r(x)$ can be deduced from the one of  $\mathcal{A}_n(x)$.
We therefore set out below the first result, which concerns the bias of $\widehat{r}_n (x)$,  for obtaining the asymptotic normality of $\widehat{r}_n (x)$.

\begin{Proposition}\label{propbr1} Under Assumptions \ref{hyph1}, \ref{hypf1} and \ref{hypk1}, we have
$$
\E_q\left[\widehat{r}_n(x)\right] - r(x)=\frac{qh_n^2 }{[2]_q f (x)} \Big(\left(D_q^2 g\right)(x) - r(x)\left(D_q^2 f\right)(x)\Big)\int_{-\nu}^{\nu} u^2K_q(u)du + O\left(q^3h_n^3b_n^{-2}\right).
$$
\end{Proposition}
\textbf{Proof.}
We first note that, for $x$ such that $f(x)$ satisfies (\ref{eqnbnf1}) and that $\widehat{f}_n (x)$ satisfies (\ref{eqnbnf2}) and  using Theorem \ref{theoconstg1} above, we get
$$
\widehat{r}_n(x) - r(x) = \left(\frac{f(x) - \widehat{f}_n (x)}{\widehat{f}_n (x)} +1 \right)\frac{\mathcal{A}_n(x)}{f (x)} = \left(o(1) +1 \right)\frac{\mathcal{A}_n(x)}{f (x)}.
$$
Now, using the following results (see \cite{badrani2024}) for $k\in \{0,1\}$
\begin{equation}\label{eqn1eqf}
\E_q\left[\widehat{\Gamma}_{k,n}(x)\right] - \Gamma_k(x)=\frac{qh_n^2 }{[2]_q} \left(D_q^2 f\right)(x)\int_{-\nu}^{\nu} u^2K_q(u)du + O\left(q^3h_n^3\right),
\end{equation}
and from (\ref{eqnbnf1}) we obtain
\begin{eqnarray*}
\E_q\left(\mathcal{A}_n(x)\right) &=& \left[\E_q\left(\widehat{g}_n (x)\right) - g(x)\right] - r(x)\left[\E_q\left(\widehat{f}_n (x)\right) - f(x)\right]\\
&=&\left[\frac{qh_n^2 }{[2]_q} \left(D_q^2 g\right)(x)\int_{-\nu}^{\nu} u^2K_q(u)du + o(qh_n^2)\right] - r(x)\left[\frac{qh_n^2 }{[2]_q} \left(D_q^2 f\right)(x)\int_{-\nu}^{\nu} u^2K_q(u)du + O\left(q^3h_n^3\right)\right]\\
&=&\frac{qh_n^2 }{[2]_q} \Big(\left(D_q^2 g\right)(x) - r(x)\left(D_q^2 f\right)(x)\Big)\int_{-\nu}^{\nu} u^2K_q(u)du + O(\|g\|_{L_q^{\infty}}q^3h_n^3b_n^{-1}).
\end{eqnarray*}
It follows that
\begin{eqnarray*}
\E_q\left[\widehat{r}_n(x)\right] - r(x)&=&  \big(o(1) +1 \big)\frac{1}{f (x)}\E_q\left(\mathcal{A}_n(x)\right)\\
&=& \frac{qh_n^2 }{[2]_q f (x)} \Big(\left(D_q^2 g\right)(x) - r(x)\left(D_q^2 f\right)(x)\Big)\int_{-\nu}^{\nu} u^2K_q(u)du + O\left(q^3h_n^3b_n^{-2}\right).\Box
\end{eqnarray*}

From this result of Proposition \ref{propbr1}, it is clear that, in view of Assumptions \ref{hyph1} and  \ref{hypf1}
$$
\E_q\left[\widehat{r}_n(x)\right] - r(x) \rightarrow 0 \, \mbox{ as }\, n \rightarrow +\infty.
$$

\medskip

In order to establish the $q$-analog of asymptotic of $\mathcal{A}_n(x)$, we need to introduce Lyapunov's central limit theorem using $q$-calculus. The standard version of Lyapunov's Theorem is found in \cite{cuzzocrea2021}. Below, we present the $q$-calculus version  proved in \cite{badrani2024}.

\begin{Theorem}\label{theoly1}
Let $0<q<1$. Considering that $\left\{X_1,\cdots,X_n\right\}$ is a sequence of independent random variables, each with finite $q$-expected value $\mu_{kq}$ and $q$-variance $\sigma_{kq}^2$. We let $V_{nq}^2=\sum_{k=1}^n \sigma_{kq}^2$. In the case when, the $q$-Lyapunov condition
\begin{equation}\label{condly1}
\lim_{n\rightarrow +\infty} \frac{1}{V_{nq}^{2+\delta}} \sum_{k=1}^n\E_q \left[\left|X_k - \mu_{kq}\right|^{2+\delta}\right]=0
\end{equation}
is satisfied for some $\delta>0$, and $q$-commutation relation:
$$
\left(X_k - \mu_{kq}\right)\left(X_\ell - \mu_{\ell q}\right)= q\left(X_\ell - \mu_{\ell q}\right)\left(X_k - \mu_{kq}\right), \, \mbox{ and }\, \sigma_{kq}^2\sigma_{\ell q}^2=q^2\sigma_{\ell q}^2\sigma_{kq}^2\, , \mbox{ for all }\, 1\leqslant k\neq \ell \leqslant n,
$$
are also satisfied, then $\frac{X_k - \mu_{kq}}{V_{nq}}$ converge in $q$-distribution to a standard normal random variable and  as $n$ goes to infinity:
$$
 \sum_{k=1}^n \frac{X_k - \mu_{kq}}{V_{nq}} \rightarrow \mathcal{N}_q(0,1),\, \mbox{ in the }\, q\mbox{-distribution sense, }
$$
where $\mathcal{N}_q(0,1)$ is the $q$-centred reduced Gaussian distribution.
\end{Theorem}

\begin{Remark}
On the $q$-commutation relation, such assumptions are often made and used in algebra and matrix calculus \citep{elkaoutit2004}. It has also been considered in \cite{badrani2024}.
\end{Remark}

\medskip
The following result states that the sequence of random variables $Z_i(x)$  satisfy the $q$-Lyapunov condition \ref{condly1}. The details of its proof are provided in Section \ref{proofs1}.

\begin{Lemma}\label{lemaz1} Under Assumptions \ref{hyph1}, \ref{hypf1}, \ref{hypy1} and \ref{hypk1}, $x\in \R$ and $0<q<1$. The sequence  $Z_1(x),\cdots, Z_n(x)$ random variables satisfies the $q$-Lyapunov condition  
$$
\lim_{n\rightarrow +\infty} \frac{1}{L_{nq}^{2+\delta}} \sum_{i=1}^n\E_q\left[\left|Z_i(x) - \E_q\big(Z_i(x)\big)\right|^{2+\delta}\right]=0
$$
for some $\delta>0$, and \, $L_{nq}^{2} = \sum_{i=1}^n Var_q\big(Z_i(x)\big)$.
\end{Lemma}

\medskip

\begin{Lemma}\label{lemaz2} Under Assumptions \ref{hyph1}, \ref{hypf1}, \ref{hypy1} and \ref{hypk1}
$$
\frac{\sum_{i=1}^n \big(Z_i(x) - \E\left(Z_i(x)\right)\big)}{\sqrt{\sum_{i=1}^n Var_q\big(Z_i(x)\big) }} \rightarrow \mathcal{N}_q(0,1) \, \mbox{ in the } q - \mbox{ distribution sense }.
$$
\end{Lemma}
\textbf{Proof.} This follows directly from Lemma \ref{lemaz1} and Theorem \ref{theoly1}.$\Box$

\bigskip

\begin{Proposition}\label{propo2an} Under Assumptions \ref{hyph1}, \ref{hypf1}, \ref{hypy1} and \ref{hypk1}
$$
\frac{\mathcal{A}_n(x)  - \E\left(\mathcal{A}_n(x) \right)}{\sqrt{Var_q\left(\mathcal{A}_n(x) \right) }} \rightarrow \mathcal{N}_q(0,1) \, \mbox{ in the } q \mbox{-} \mbox{distribution sense }.
$$
\end{Proposition}
\textbf{Proof.} Consider  writing of $\mathcal{A}_n(x)$ in (\ref{eqn1anz}). The proof of the proposition comes down to proving that
$$
\frac{\sum_{i=1}^n \big(Z_i(x) - \E\left(Z_i(x)\right)\big)}{\sqrt{\sum_{i=1}^n Var_q\big(Z_i(x)\big) }} \rightarrow \mathcal{N}_q(0,1) \, \mbox{ in the } q\mbox{-distribution sense }.
$$
Therefore, Lemma \ref{lemaz2} completes the demonstration. $\Box$

In the following, we define the sequence
$$
\mathcal{B}_n(x) = \frac{f(x)\sqrt{nh_n}\left(\widehat{r}_n (x) - r(x)\right)}{\sqrt{Var_q \left(Y | X=x\right)f(x)\int_{-\nu}^{+\nu} K_q^2(u)d_qu }}.
$$
From (\ref{eqn1lq}) we have
$$
\mathcal{B}_n(x) = \frac{f(x)}{\widehat{f}_n (x)}\frac{\mathcal{A}_n(x)}{\sqrt{\sum_{i=1}^n Var_q\big(Z_i(x)\big)}} = \left(\frac{f(x) - \widehat{f}_n (x) }{\widehat{f}_n (x)}+1\right)\frac{\mathcal{A}_n(x)}{\sqrt{\sum_{i=1}^n Var_q\big(Z_i(x)\big)}}.
$$
Now, using Theorem \ref{theoconstg1} above, Assumption \ref{hyph1} and (\ref{eqnbnf2}), it follows that
\begin{equation}\label{eq3bn}
\mathcal{B}_n(x) =\big(o(1)+1\big)\frac{\mathcal{A}_n(x)}{\sqrt{ Var_q\big(\mathcal{A}_n(x)\big)}}.
\end{equation}
This leads to the following lemma:

\begin{Lemma}\label{lemaB1} Under Assumptions \ref{hyph1}, \ref{hypf1}, \ref{hypy1} and \ref{hypk1}
$$
\frac{\mathcal{B}_n(x)  - \E\left(\mathcal{B}_n(x) \right)}{\sqrt{Var_q\left(\mathcal{B}_n(x) \right) }} \stackrel{\mathscr{L}}{\longrightarrow } \mathcal{N}_q(0,1) \,\mbox{ as }\, n\rightarrow +\infty \mbox{ in the } q\mbox{-distribution sense }.
$$
\end{Lemma}
In the sequel, we present the results establishing the almost sure convergence of the estimators $\widehat{f}_n(x)$, $\widehat{g}_n(x)$, and $\widehat{r}_n(x)$. Additionally, we determine their convergence rates. Proofs not provided in this section, as the proof of below lemma, they are presented in Section \ref{proofs1}.

\begin{Lemma}\label{lemmaconst2} Under Assumptions  \ref{hyph1}, \ref{hypy1} and  \ref{hypk1}, for all $t>0$ and $n\geqslant 1$, for $k\in\{0,1\}$, we have
$$
\sup_{x\in \R}P_q\left\{\left|\widehat{\Gamma}_{k,n}(x)-\E_q\left(\widehat{\Gamma}_{k,n}(x)\right)\right| \geqslant t\right\}\leqslant \exp_q\left\{-\frac{t^2}{[2]_q(v_{(n,k)}+ct)}\right\},
$$
where, $c=\frac{M\left(\left(c_0  \ln_q (n)\right)^{\frac{k}{2}} \|f\|_{L_q^{\infty}} + o(1)\right)[2]_q}{[3]_q! nh_n}$, $v_{(n,k)}$ is given in (\ref{eqn1vnk}) and the constant $M$ is given in  Assumptions  \ref{hypk1}.
\end{Lemma}

\medskip

\begin{Proposition}\label{propqg1} Under Assumptions  \ref{hyph1}, \ref{hypy1} and  \ref{hypk1}, for all $t>0$ and $n\geqslant 1$ and for $k\in\{0,1\}$, we have
$$
\sup_{x\in\R}\left|\widehat{\Gamma}_{k,n}(x)-\E_q\left(\widehat{\Gamma}_{k,n}(x)\right)\right|= O_{a.s}\left(\sqrt{v_{(n,k)}\ln_q n}\right).
$$
where $v_{(n,k)}$ is given in (\ref{eqn1vnk}).
\end{Proposition}
\textbf{Proof.} In the following of Lemma \ref{lemmaconst2},
\begin{eqnarray*}
\frac{v_{(n,k)}}{c}&=&\frac{\left(\left(c_0  \ln_q (n)\right)^{\frac{k}{2}}\|f\|_{L_q^{\infty}} + o(1)\right)}{nh_n} \int_{-\nu}^{\nu} K_q^2(u)d_qu \times \frac{[3]_q! nh_n}{M\left(\left(c_0  \ln_q (n)\right)^{\frac{k}{2}}\|f\|_{L_q^{\infty}} + o(1)\right)[2]_q}\\
& = &\frac{[3]_q!}{M[2]_q} \int_{-\nu}^{\nu} K_q^2(u)d_qu=c_1,
\end{eqnarray*}
we get
$$
\exp_q\left\{-\frac{t^2}{[2]_q(v+ct)}\right\} = \exp_q\left\{-\frac{t^2}{[2]_qc(c_1+t)}\right\}
$$
Taking $t:=t_n = L\sqrt{v_{(n,k)}\ln_q n}$ where $L$ is a constant such that $L>\sqrt{2[2]_q}$. With this choice, we have $t_n \rightarrow 0$ as $n\rightarrow +\infty$, for $n$ large enough such that $t_n \leqslant c_1$ we get
$$
\exp_q\left\{-\frac{t_n^2}{[2]_qc(c_1+t_n)}\right\} = \exp_q\left\{-\frac{ L^2 v_{(n,k)}\ln_q n}{[2]_qc(c_1+t_n)}\right\} \leqslant \exp_q\left\{-\left[\frac{ L^2}{2[2]_q}\right] \ln_q n\right\}.
$$																									
$[x]$ stands  the integer part of $x$. Now, in view of the third equality of (\ref{eqnlnq3})	and Equation (25) in \cite{tsallis1994}, the above inequality becomes    
$$
\exp_q\left\{-\frac{t_n^2}{[2]_qc(c_1+t_n)}\right\} \leqslant\left(\exp_q\left\{ \ln_q n\right\}\right)^{-\left[\frac{ L^2}{2[2]_q}\right]} \leqslant n^{-\left[\frac{ L^2}{2[2]_q}\right]}.
$$	
Then, thanks to the Borel Cantelli lemma,the proof is completed. $\Box$

\bigskip

\section{Proofs}\label{proofs1}
In this last section, we complete with the proofs of the stated results. We begin with some necessary lemmas.

\subsection{Technical lemmas}

\begin{Lemma} Under Assumptions \ref{hypf1} and \ref{hypk1},  we have
$$
\E_q\big(Z_i(x)\big)=\frac{qh_n^2 }{[2]_q n} \big(\left(D_q^2 g\right)(x) - r(x)\left(D_q^2 f\right)(x)\big)\int_{-\nu}^{\nu} u^2K_q(u)du + O\left(\frac{q^3h_n^3}{nb_n}\right).
$$
\end{Lemma}
\textbf{Proof.}
We have 
\begin{eqnarray*}
\E_q\left(\big(Y_i - r(x)\big) K_q\left(\frac{x-X_i}{h_n}\right)\right)  &=&  \int_{\R}\int_{x-h_n\nu}^{x+h_n\nu} \big(v - r(x)\big) K_q\left(\frac{x-u}{h_n}\right) f_{_{X,Y}}(u,v)d_qu \, d_qv\\
&=&  h_n\int_{-\nu}^{\nu} g(x- h_nu) K_q(u)d_qu - r(x)h_n \int_{-\nu}^{\nu}K_q(u)  f(x- h_nu )d_qu.
\end{eqnarray*}
Applying a second-order $q$-Taylor series expansion about $x$ results in $g(x- h_nu)$, we obtain
$$
\int_{-\nu}^{\nu} g(x- h_nu) K_q(u)d_qu = g(x)+ \frac{qh_n^2 }{[2]_q} \left(D_q^2 g\right)(x)\int_{-\nu}^{\nu} u^2K_q(u)d_qu + O(q^3h_n^3)
$$
and
$$
\int_{-\nu}^{\nu} f(x- h_nu) K_q(u)d_qu = f(x) + \frac{qh_n^2 }{[2]_q} \left(D_q^2 f\right)(x)\int_{-\nu}^{\nu} u^2K_q(u)d_qu + O(q^3h_n^3),
$$
it follows that, using (\ref{eqnbnf1}) and the fact that $g$ is bounded
$$
\E_q\left(\big(Y_i - r(x)\big) K_q\left(\frac{x-X_i}{h_n}\right)\right) = \frac{qh_n^3 }{[2]_q} \big(\left(D_q^2 g\right)(x) - r(x)\left(D_q^2 f\right)(x)\big)\int_{-\nu}^{\nu} u^2K_q(u)du + O\left(q^3b_n^{-1}h_n^4\right).
$$
And the result is obtained.$\Box$

\medskip

\begin{Lemma}\label{lemavarq1} Under Assumptions \ref{hypf1} and \ref{hypk1},  we have
$$
Var_q\big(Z_i(x)\big)=\frac{1}{n^2 h_n} Var_q \left(Y | X=x\right)f(x)\int_{-\nu}^{+\nu} K_q^2(u)d_qu + O\left(\frac{qh_n}{n^2}\right).
$$
\end{Lemma}
\textbf{Proof.} We have for $k\in\N$
\begin{equation}\label{qnro1}
\E_q\left[Y_i^k K_q^2\left(\frac{x-X_i}{h_n}\right)\right]=h_n \int_{-\nu}^{+\nu} K_q^2(u) \Gamma_k(x- h_nu) d_qu.
\end{equation}
Since, the $q$-Taylor formula of $\Gamma_2$ gives
$$
\Gamma_k(x- h_nu) = \Gamma_k(x)-h_n u\left(D_q^1 \Gamma_k \right)(x) + \frac{q h_n^2 u^2}{[2]_q} \left(D_q^2 \Gamma_k\right)(x) + O(q^3h_n^3u^3), 
$$
it follows that, using the fact that $\int_{-\nu}^{+\nu}u K^2(u) d_qu =0$,
$$
\E_q\left[Y_i^k K^2\left(\frac{x-X_i}{h_n}\right)\right]=h_n \Gamma_k(x)  \int_{-\nu}^{+\nu} K_q^2(u)d_qu  + O\left(qh_n^3\right).
$$
$\Gamma_k(x)$ was defined in (\ref{eqnrho1}).

It follows that
\begin{eqnarray*}
& &\E_q\left[\big(Y_i - r(x)\big)^2 K_q^2\left(\frac{x-X_i}{h_n}\right)\right]\\
& & \quad =  \E_q\left[Y_i^2K^2\left(\frac{x-X_i}{h_n}\right)\right]+r(x)^2\E_q\left[K_q^2\left(\frac{x-X_i}{h_n}\right)\right]- 2r(x)\E_q\left[Y_iK_q^2\left(\frac{x-X_i}{h_n}\right)\right]\\
& &  \quad = h_n \Gamma_2(x)  \int_{-\nu}^{+\nu} K_q^2(u)d_qu + r(x)^2h_n \Gamma_0(x)  \int_{-\nu}^{+\nu} K_q^2(u)d_qu -2r(x)  h_n \Gamma_1(x)  \int_{-\nu}^{+\nu} K_q^2(u)d_qu+ O\left(qh_n^3\right)\\
& &  \quad = h_n\left(\Gamma_2(x) + r(x)^2 f(x)-2r(x)g(x)\right)\int_{-\nu}^{+\nu} K_q^2(u)d_qu + O\left(qh_n^3\right)\\
& &  \quad = h_n\big(\Gamma_2(x) + r(x)^2 f(x)-2r(x)^2 f(x)\big)\int_{-\nu}^{+\nu} K_q^2(u)d_qu + O\left(qh_n^3\right)\\
& &  \quad = h_n\big(\Gamma_2(x) - r(x)^2 f(x)\big)\int_{-\nu}^{+\nu} K_q^2(u)d_qu + O\left(qh_n^3\right)\\
& &  \quad = h_n\big(\E_q\left(Y^2 | X=x\right)  - \E_q^2\left(Y | X=x\right)\big)f(x)\int_{-\nu}^{+\nu} K_q^2(u)d_qu + O\left(qh_n^3\right)\\
& &  \quad = h_n Var_q \left(Y | X=x\right)f(x)\int_{-\nu}^{+\nu} K_q^2(u)d_qu + O\left(qh_n^3\right).\\
\end{eqnarray*}
Now,
\begin{eqnarray*}
& &Var_q\left[\big(Y_i - r(x)\big) K_q\left(\frac{x-X_i}{h_n}\right)\right]\\
& & \quad = \E_q\left[\big(Y_i - r(x)\big)^2 K_q^2\left(\frac{x-X_i}{h_n}\right)\right] - \left(\E_q\left[\big(Y_i - r(x)\big) K_q\left(\frac{x-X_i}{h_n}\right)\right]\right)^2 \\
& & \quad =  h_n Var_q \left(Y | X=x\right)f(x)\int_{-\nu}^{+\nu} K_q^2(u)d_qu\\
& & \quad \quad - \left(\frac{qh_n^3 }{[2]_q} \big(\left(D_q^2 g\right)(x) - r(x)\left(D_q^2 f\right)(x)\big)\int_{-\nu}^{\nu} u^2K_q(u)du\right)^2 + O\left(h_n^6\right)+ o\left(h_n^6b_n^{-2}\right) + O\left(qh_n^3\right)\\
& & \quad =  h_n Var_q \left(Y | X=x\right)f(x)\int_{-\nu}^{+\nu} K_q^2(u)d_qu + O\left(qh_n^3\right)
\end{eqnarray*}
Hence, the proof is complete. $\Box$

\medskip

\subsection{Proof of Lemma \ref{lemaz1}}

Setting $\delta=1$, we have $\E\left[\left|Z_i(x) - \E_q\big(Z_i(x)\big)\right|^3\right] \leqslant   8\E\left[\left|Z_i^3(x)\right|\right]$.
However, getting $\E\left[\left|Z_i^3(x)\right|\right]$ is getting $\E\left[\left|\big(Y_i - r(x)\big) K_q\left(\frac{x-X_i}{h_n}\right)\right|^3\right]$.\\
In view of assumptions, one can see the following inequality
$$
\E\left[\left|\big(Y_i - r(x)\big) K_q\left(\frac{x-X_i}{h_n}\right)\right|^3\right] \leqslant \left( \max_{1\leqslant i \leqslant n}|Y_i| + \sup_{x\in \R}|r(x)|\right) h_n\int_{-\nu}^{+\nu} K_q^3(u)f(x-h_n u)d_qu.
$$
Since
$$
h_n \int_{-\nu}^{+\nu} K_q^3(u)f(x-h_n u)d_qu=h_n f(x)\int_{-\nu}^{+\nu} K_q^3(u)d_qu + O\left(qh_n^3\right),
$$
it follows from Assumptions \ref{hyph1}, \ref{hypf1}, \ref{hypy1} and \ref{hypk1},
$$
\E\left[\left|\big(Y_i - r(x)\big) K_q\left(\frac{x-X_i}{h_n}\right)\right|^3\right] \leqslant h_n\left(A_{nq} + b_n^{-1}\|g\|_{L_q^{\infty}}\right).
$$
Then
$$
\E\left[\left|Z_i(x) - \E_q\big(Z_i(x)\big)\right|^3\right] \leqslant \frac{1}{n^3h_n^2}\left(A_{nq} + b_n^{-1}\|g\|_{L_q^{\infty}}\right),
$$
and
$$
\sum_{i=1}^n \E\left[\left|Z_i(x) - \E_q\big(Z_i(x)\big)\right|^3\right] \leqslant \frac{1}{n^2h_n^2}\left(A_{nq} + b_n^{-1}\|g\|_{L_q^{\infty}}\right).
$$
Now, from Lemma \ref{lemavarq1},
\begin{equation}\label{eqn1lq}
L_{nq}^{2} = \frac{1}{n h_n} Var_q \left(Y | X=x\right)f(x)\int_{-\nu}^{+\nu} K_q^2(u)d_qu + O\left(\frac{qh_n}{n}\right),
\end{equation}
we have
$$
L_{nq}^{3} = \left(\frac{1}{n h_n} Var_q \left(Y | X=x\right)f(x)\int_{-\nu}^{+\nu} K_q^2(u)d_qu + O\left(\frac{qh_n}{n}\right)\right)^{\frac{3}{2}}.
$$
And, from Assumptions \ref{hyph1} and \ref{hypy1}, we have $\frac{A_{nq} + b_n^{-1}\|g\|_{L_q^{\infty}}}{\sqrt{nh_n}} \rightarrow 0$ as $n \rightarrow + \infty$. Hence, the proof of this lemma is completed.$\Box$

\medskip

\begin{Lemma}[$q$-Markov's Inequality]
Let $X$ a random variable with a $q$-distribution. If $X$ is a nonnegative random variable and $a > 0$, then we have
$$
P_q\left(X\geqslant a \right)\leqslant \frac{\E_q(X)}{a}.
$$
\end{Lemma}
\textbf{Proof.} Clearly, we have $\E_q(X) \geqslant a P_q\left(X\geqslant a \right)$.$\Box$

\medskip

The lemma below requires two conditions of $q$-commutation, conditions necessary for $q$-calculus operations. These conditions have been considered in \cite{badrani2024}.

\begin{Lemma}[$q$-Bernstein's Inequality]\label{lembersntein1}
Let $0 < q < 1$. Considering  $X_1; \cdots ;X_n$ be $n$ independent and identically variables $q$-distributed of  random
variables of $X$ with finite $q$-expected value $\E_q(X_i)$.  Let $Z=\sum_{i=1}^n X_i$. We assume that the $q$-commutation relations $\E_q(X_i)\E_q(X_j) = q\E_q(X_j)\E_q(X_i)$ and $\E_q(X_i^2)\E_q(X_j^2) = q^2\E_q(X_j^2)\E_q(X_i^2)$ are satisfied and that there exist two constants $v>0$ and $c>0$ such that $\sum_{i=1}^n \E_q\left(X_i^2\right) \leqslant v$ and 
$$
\forall k\geqslant 3,\,\sum_{i=1}^n \E_q\left[\left(X_i\right)_+^k\right] \leqslant \frac{v [k]_q!c^{k-2}}{[2]_q}. 
$$
Then for all $\lambda \in [0,1/c[$,
$$
\ln_q \E_q\left[\exp_q\left\{\lambda\big(Z-\E_q (Z)\big)\right\}\right] \leqslant \frac{v\lambda^2}{[2]_q(1-c\lambda)}.
$$
which implies that for all $t\geqslant 0$,
\begin{equation}\label{eqnq1ber}
P_q\left\{\left|Z-\E_q(Z)\right| \geqslant t\right\}\leqslant \exp\left\{-\frac{t^2}{[2]_q(v+ct)}\right\}.
\end{equation}
\end{Lemma}
\textbf{Proof.} We have
$\lambda\left(Z-\E_q(Z)\right)=\sum_{i=1}^n\left(\lambda X_i - \E_q(\lambda X_i)\right)$. Using (\ref{eqnlnq2}) we obtain
$$
\exp_q \left\{\lambda\left(Z-\E_q(Z)\right)\right\}=\exp_q\left\{\sum_i^n  \big(\lambda X_i - \E_q(\lambda X_i)\big)\right\} = \prod_i^n \exp_q\left\{\lambda X_i - \E_q(\lambda X_i)\right\}.
$$
It follows that
$$
\ln_q\Big(\E_q\left[\exp_q\left\{\lambda\left(Z-\E_q(Z)\right)\right\}\right]\Big)=\sum_{i=1}^n\left\{\ln_q\Big(1+ \lambda \E_q(X_i) +\E_q\left[\phi(\lambda X_i)\right] \Big)  - \lambda \E_q(X_i)\right\},
$$
where $\phi(u)=\exp_q(u) -u-1$.\\
Now, by considering (\ref{eqnlnq1}), we obtain
$$
\ln_q\Big(\E_q\left[\exp\left\{\lambda\left(Z-\E_q(Z)\right)\right\}\right]\Big) \leqslant \sum_{i=1}^n \E_q\left[\phi(\lambda X_i)\right].
$$
As for all $u\geqslant 0$, $\phi(u)\leqslant \frac{u^2}{[2]_q}$ and that $\phi(0)=0$, we have
$$
\phi(\lambda X_i) = \phi((\lambda X_i)_{-}) + \phi((\lambda X_i)_{+}) \leqslant \frac{\lambda^2 (X_i)^2_{-}}{[2]_q} + \sum_{k\geqslant 2} \frac{\lambda^k(X_i)^k_{+}}{k!} = \frac{\lambda^2 (X_i)^2}{[2]_q} + \sum_{k\geqslant 3} \frac{\lambda^k(X_i)^k_{+}}{[k]_q!}.
$$
So, for $\lambda \in[0,\, 1/c[$,
$$
\E_q\left[\phi(\lambda X_i)\right] \leqslant \frac{\lambda^2 \E_q\left[(X_i)^2\right]}{[2]_q} + \sum_{k\geqslant 3} \frac{\lambda^k\E_q\left[(X_i)^k_{+}\right]}{[k]_q!},
$$
consequently 
\begin{eqnarray*}
\sum_{i=1}^n \E_q\left[\phi(\lambda X_i)\right] &\leqslant& \frac{\lambda^2}{[2]_q}\sum_{i=1}^n \E_q\left[(X_i)^2\right] + \sum_{k\geqslant 3}\frac{\lambda^k}{[k]_q!} \sum_{i=1}^n  \E_q\left[(X_i)^k_{+}\right]\\ 
 &\leqslant& \frac{\lambda^2 v}{[2]_q} + \sum_{k\geqslant 3} \frac{v \lambda^k c^{k-2}}{[2]_q}\\
&\leqslant& \frac{\lambda^2 v}{[2]_q(1-c\lambda)}.
\end{eqnarray*}
It follows
$$
\ln_q\Big(\E_q\left[\exp\left\{\lambda\left(Z-\E_q(Z)\right)\right\}\right]\Big) \leqslant \frac{\lambda^2 v}{[2]_q(1-c\lambda)},
$$
that implies
$$
\E_q\left[\exp\left\{\lambda\left(Z-\E_q(Z)\right)\right\}\right] \leqslant \exp_q\left\{\frac{\lambda^2 v}{[2]_q(1-c\lambda)}\right\}.
$$
Now, according to $q$-Markov's Inequality, we get
\begin{eqnarray}
P_q\left(Z-\E_q(Z) \geqslant t \right) &=& P_q\big(\E_q\left[\exp_q \left\{\lambda\left(Z-\E_q(Z)\right)\right\}\right]\nonumber\\
& \geqslant& \exp_q \left\{\lambda t\right\}\big)\leqslant \exp_q\left\{-\lambda t + \frac{v\lambda^2}{[2]_q(1-c\lambda)}\right\}.\label{eqnq3markov}
\end{eqnarray}
The optimal choice of $\lambda$ to minimize $\exp_q\left\{-\lambda t + \frac{v\lambda^2}{[2]_q(1-c\lambda)}\right\}$ is gives for \,$\lambda = \frac{t}{[2]_q(v+ct)}$.\\
By substituting $\lambda$ in (\ref{eqnq3markov}), we obtain (\ref{eqnq1ber}) after simplification.$\Box$

\subsection{Proof of Lemma \ref{lemmaconst2}}
Let $Z=\frac{1}{nh_n}\sum_{i=1}^n Y_iK_q\left(\frac{x-X_i}{h_n}\right) = \sum_{i=1}^n T_i$, with $T_i=\frac{1}{nh_n}Y_iK_q\left(\frac{x-X_i}{h_n}\right)$.\\
We have from (\ref{qnro1}),
$$
\E_q\left[Y^2K_q^2\left(\frac{x-X}{h_n}\right)\right] = h_n\int_{-\nu}^{\nu} K_q^2(u)\left(\Gamma_2(x) + o(1)\right)d_qu.
$$
Then
$$
\sum_{i=1}^n\E_q\left(T_i^2\right)=\frac{1}{nh_n} \int_{-\nu}^{\nu} K_q^2(u)\left(\Gamma_2(x) + o(1)\right)d_qu =v.
$$
On the other hand, for $k\geqslant 3$,
$$
\E_q\left[Y^kK^k_q\left(\frac{x-X_i}{h_n}\right)\right] = \int_{-\nu}^{\nu} K_q^k(u)\Gamma_k (x-h_nu)d_qu.
$$
On the other hand, first-order $q$-Taylor series expansion of $\Gamma_k (x-h_nu)$ gives
$$
\Gamma_k (x-h_nu) = \Gamma_k (x)-h_n u \left(D_q^1 \Gamma_k \right)(x)+ o(h_n),
$$
then
$$
\E_q\left[K^k_q\left(\frac{x-X_i}{h_n}\right)\right] = h_n\int_{-\nu}^{\nu} K_q^k(u)\left( \Gamma_k(x) + o(1)\right)d_qu.
$$
It follows that, using Assumption \ref{hypk1}, we get
\begin{eqnarray*}
\sum_{i=1}^n\E_q\left(T_i^k\right)&=&\frac{1}{n^{k-1}h_n^{k-1}} \int_{-\nu}^{\nu} K_q^k(u)d_qu\\ 
                                  &=&\frac{\left(\Gamma_k(x) + o(1)\right)}{n^{k-1}h_n^{k-1}} \int_{-\nu}^{\nu} K_q^2(u) K_q^{k-2}(u)d_qu\\
																	&\leqslant &\left(\frac{\left(\Gamma_k(x) + o(1)\right)[2]_q}{[k]_q! }\right) \left(\frac{M}{nh_n}\right)^{k-2} \frac{ [k]_q! v}{[2]_q}\\
																	&\leqslant &\left(\left(\frac{\left(\Gamma_k(x) + o(1)\right)[2]_q}{[k]_q! }\right)^{\frac{1}{k-2}}\frac{M}{nh_n} \right)^{k-2} \frac{ [k]_q! v}{[2]_q}\\
																	&\leqslant &\left(\frac{\left(\Gamma_k(x) + o(1)\right)[2]_q}{[3]_q! }\frac{M}{nh_n} \right)^{k-2} \frac{ [k]_q! v}{[2]_q}.
\end{eqnarray*}
Then
$$
\sum_{i=1}^n\E_q\left(T_i^k\right) \leqslant \frac{c^{k-2} [k]_q! v}{[2]_q}.
$$ 
Therefore, the Proposition is proved with Lemma \ref{lembersntein1}. $\Box$

\subsection{Proof of Theorem \ref{Theonr1}}

From Lemma \ref{lemaB1}, for $n$ large enough \, $\mathcal{B}_n(x)$ follows the normal distribution \, $\mathcal{N}_q\Big(\E_q\left(\mathcal{B}_n(x)\right)\, ,\, Var_q\left(\mathcal{B}_n(x)\right)\Big)$. Then, setting 
$$
\gamma(x)= f(x)^{1/2}\left(Var_q \left(Y | X=x\right)\int_{-\nu}^{+\nu} K_q^2(u)d_qu\right)^{-1/2},
$$
for $n$ large enough
$$
\sqrt{nh_n}\left(\widehat{r}_n (x) - r(x)\right)  \mbox{  follows the normal distribution } \mathcal{N}_q\left(\frac{\E_q\left(\mathcal{B}_n(x)\right)}{\gamma(x)}\, ,\, \frac{Var_q\left(\mathcal{B}_n(x)\right)}{\gamma(x)^2}\right).
$$
However, using Proposition \ref{propbr1}
\begin{eqnarray*}
\E_q\left(\mathcal{B}_n(x)\right) &=& \gamma(x)\sqrt{nh_n}\left(\E_q\left(\widehat{r}_n (x)\right) - r(x)\right)\\
&=& \gamma(x) \frac{q\sqrt{nh_n^5} }{[2]_q f (x)} \Big(\left(D_q^2 g\right)(x) - r(x)\left(D_q^2 f\right)(x)\Big)\int_{-\nu}^{\nu} u^2K_q(u)du + O\left(q \sqrt{nh_n^5}h_n b_n^{-2}\right),
\end{eqnarray*}
and in fact that $\sqrt{nh_n^5}\rightarrow c$, it follows that
$$
\E_q\left(\mathcal{B}_n(x)\right) =  \frac{q c \,\gamma(x)}{[2]_q f (x)} \Big(\left(D_q^2 g\right)(x) - r(x)\left(D_q^2 f\right)(x)\Big)\int_{-\nu}^{\nu} u^2K(u)du + O\left(q c h_n b_n^{-2}\right).
$$
Now, from (\ref{eq3bn}), We have equality $Var_q\left(\mathcal{B}_n(x)\right)=1$. Then, in view of Assumption \ref{hyph1}, when $n\rightarrow +\infty$, we have
\begin{multline*}
\mathcal{N}_q\left(\frac{\E_q\left(\mathcal{B}_n(x)\right)}{\gamma(x)}\, ,\, \frac{Var_q\left(\mathcal{B}_n(x)\right)}{\gamma(x)^2}\right)\\
= \mathcal{N}_q\left(\frac{q c}{[2]_q f (x)} \Big(\left(D_q^2 g\right)(x) - r(x)\left(D_q^2 f\right)(x)\Big)\int_{-\nu}^{+\nu} u^2K_q(u)d_qu \, , \, f(x)^{-1}\left(Var_q \left(Y | X=x\right)\int_{-\nu}^{+\nu} K_q^2(u)d_qu\right)\right).\Box
\end{multline*}

\subsection{Proof of Theorem \ref{theoconstr1}}
From (\ref{eqnbnf1}), we have
\begin{eqnarray*}
\left|\widehat{r}_n(x)-r(x)\right| &=& \frac{\left|\left(\widehat{g}_n(x)-g(x)\right)f(x) + \left(\widehat{f}_n(x)-f(x)\right)g(x)\right|}{f(x)\widehat{f}_n(x)}\\
                               &\leqslant & b_n^{-1}\frac{\left|\left(\widehat{g}_n(x)-g(x)\right)f(x) + \left(\widehat{f}_n(x)-f(x)\right)g(x)\right|}{\widehat{f}_n(x)}\\
\end{eqnarray*}
from Assumption \ref{hyph1}(ii), we have  $\left|\widehat{f}_n (x) - f(x)\right| = o(b_n)$. Therefore, for $n$ large enough  we get $\left|\widehat{f}_n (x) - f(x)\right| \leqslant b_n/2$. So, $b_n \leqslant f(x) \leqslant  \left|\widehat{f}_n (x) - f(x)\right| + \widehat{f}_n (x) \leqslant b_n/2+ \widehat{f}_n (x)$.\\
 Hence, 
\begin{equation}\label{eqnbnf2}
\widehat{f}_n (x) \geqslant b_n/2.
\end{equation}
 It follows that
$$
\left|\widehat{r}_n (x) - r(x)\right| \leqslant 2 b_n^{-2} \|f\|_{L_q^{\infty}}\left|\widehat{g}_n (x) - g(x)\right|  +  b_n^{-2}\|g\|_{L_q^{\infty}}\left|\widehat{f}_n (x) - f(x)\right|.
$$ 
Consequently,
$$
\sup_{x\in\R^p}\left|\widehat{r}_n (x) - r(x)\right| \leqslant 2 e_n^{-2} \left(\|f\|_{L_q^{\infty}} \sup_{x\in\R^p}\left|\widehat{g}_n (x) - g(x)\right|  +  \|g\|_{L_q^{\infty}} \sup_{x\in\R^p}\left|\widehat{f}_n (x) - f(x)\right|\right),
$$ 
and the proof is completed by using Theorem \ref{theoconstg1}. $\Box$

\end{document}